\magnification=1100

\def\R{\hbox{\bf R}}
\def\Z{\hbox{\bf Z}}
\def\N{\hbox{\bf N}}
\def\C{\hbox{\bf C}}
\def\P{\hbox{\bf P}}

\newdimen\sumdim

\def\diagram#1\enddiagram
{\vcenter{\offinterlineskip
\def\tvi{\vrule height 10pt depth 10pt width 0pt}
\halign{&\tvi\kern 5pt\hfil$\displaystyle##$\hfil\kern 5pt
\crcr #1\crcr}}}

\def\hfl[#1][#2][#3]#4#5{\kern-#1
\sumdim=#2 \advance\sumdim by #1 \advance\sumdim by #3
\smash{\mathop{\hbox to \sumdim{\rightarrowfill}}
\limits^{\displaystyle#4}_{\displaystyle#5}}
\kern-#3}

\def\vfl[#1][#2][#3]#4#5%
{\sumdim=#1 \advance\sumdim by #1 \advance\sumdim by #3
\smash{\mathop{\hbox to \sumdim{\rightarrowfill}}
\setbox1=\hbox{$\left\downarrow\vbox to .5\sumdim{}\right.$}
\setbox1=\hbox{\llap{$\displaystyle#4$}}
\rlap{$\displaystyle#5$}
\vcenter{\kern-#1\box1\kern-#3}}

\def\vspace[]{\noalign\vskip #1}}

\centerline{\bf Dual elliptic structures on $\C\P^2$}
\medskip\centerline{Jean-Claude Sikorav}
\bigskip
\noindent{\bf Abstract.}\hskip3mm We define the notion of tame elliptic structure $E$
on $\C\P^2$, which 
generalizes an almost complex structure $J$ on which the standard symplectic form is positive. 
An
$E$-curve is a surface in $V$ which is everywhere tangent to $E$ 
($J$-holomorphic in the almost complex case), and an
$E$-line is an $E$-curve
of degree $1$. We show that
the space $V^*$ of $E$-lines is itself a complex projective plane with a tame elliptic structure $E^*$. Moreover, 
to each $E$-curve one can associate
its dual in $V^*$, which is an $E^*$-curve.
This implies that the $E$-curves, and in particular the $J$-curves, satisfy the Pl\"ucker formulas,
which restricts their possible sets of singularities.

\bigskip
\noindent AMS 32Q65, 53C15, 53C42, 53D35, 57R17, 58J60
\vskip10mm
\noindent{\bf Introduction}
\bigskip
Let $V$ be a smooth oriented $4$-manifold, which is a rational homology $\C\P^2$ (ie $b_2(V)=1$), and let 
$J$ be an 
almost complex structure on $V$ 
which is homologically equivalent to the standard
structure $J_0$ on $\C\P^2$. This means that there is an isomorphism $H^*(V)\to H^*(\C\P^2)$ 
(rational coefficients) which is positive on $H^4$ and sends the Chern class $c_1(J)$ to $c_1(J_0)$.

\smallskip
By definition, a {\it $J$-line} is a 
$J$-holomorphic curve ($J$-curve for short) of degree $1$. By 
the {\it positivity of intersections} [McD2], it is an embedded sphere. We denote by $V^*$ the set of $J$-lines.
\smallskip
Now assume that $J$ is {\it tame}, ie positive with respect to some symplectic form $\omega$, and also that
$V^*$ is  nonempty.
Then M. Gromov 
[G] (2.4.A) [cf. also [McD1]) has proved that by two distinct
points $x,y\in V$ there passes a unique $J$-line $L_{x,y}\in V^*$, depending smoothly on $(x,y)$; also, for any given $P\in Gr_1^J(TV)$, the Grassmannian of
$J$-complex lines in $TV$, there exists a unique $J$-line $L_P\in V^*$
tangent to $P$. Furthermore, 
$V$ is oriented diffeomorphic to $\C\P^2$, $\omega$ is isomorphic to $\lambda\omega_0$ for some
 positive $\lambda$ so that $J$ is homotopic to $J_0$. Finally, $V^*$ has a natural structure of compact oriented 
$4$-manifold; although it is not explicitly stated in [G], the above properties of $V^*$ imply 
that it is also oriented diffeomorphic to 
$\C\P^2$.

\smallskip
Later, Taubes [T] proved 
that the hypothesis that $V^*$ be nonempty is unnecessary, so that all the above results hold when $J$ is tame. We shall call 
$(V,J)$ a {\it tame almost complex projective plane}.
\medskip
Following [G, 2.4.E], these facts can be extended to the case of an {\it elliptic structure} on $V$, ie one replaces
$Gr_1^J(TV)$ by a
suitable submanifold $E$ of the Grassmannian of oriented $2$-planes $Gr_2TV$. Such a structure is associated to a {\it twisted
almost complex structure} $J$, which is a fibered map from $TV$ to itself satisfying $J_v^2=-{\rm Id}$ 
but such that $J_v$ is not
necessarily linear.
\smallskip
An elliptic  structure on $V$ gives rise to a 
notion of {\it $E$-curve}, ie a surface $S\subset V$
(not necessarily embedded or immersed)
whose tangent plane at every point is an element of $E$ (for the precise definitions, see section 2). 
It will be called {\it tame} it there exists
a symplectic form $\omega$ strictly positive on each $P\in E$.
\medskip
In Gromov's words, ``all facts on $J$-curves [proved
in [G]] extend to $E$-curves with an obvious change of terminology''. In particular, 
let $V$ be a rational homology 
$\C\P^2$ equipped with a tame elliptic structure $E$ so that $(V,E)$ is homologically equivalent
to $(\C\P^2,Gr_1^{\C}(T\C\P^2))$. Then one can define the space $V^*$ of {$E$-lines} ($E$-curves of degree $1$), 
and prove that all the above properties still hold (see section 3). In particular, $V$ and $V^*$ are  
oriented diffeomorphic to
$\C\P^2$. 
\medskip
We shall call $(V,E)$ with the above properties a {\it tame elliptic projective plane}.  If 
$C\subset V$ is an $E$-curve, we define its {\it dual} 
$C^*\subset V^*$ by $C^*=\{L_{T_vC}\mid v\in C\}$. A more precise definition is given in section 4; note that one must require
that no component of $C$ is contained in an $E$-line. The main new result 
of this paper is then the following.
\medskip

\noindent{\bf Theorem.} \hskip3mm {\it Let $(V,E)$ be a tame elliptic projective plane.
Then there exists a unique elliptic structure $E^*\subset Gr_2TV^*$ on $V^*$ with the following property: 
if $C\subset V$ is an $E$-curve, then its dual $C^*\subset V^*$ is an $E^*$-curve.
\smallskip

Furthermore, $(V^*,E^*)$ is again a tame elliptic projective plane. 
Finally, the bidual  $(V^{**},E^{**})$ can be  canonically identified with $(V,E)$, and $C^{**}=C$ for every
$E$-curve $S$.}
\bigskip
If $E$ comes from an almost complex structure, one may wonder if this is also the case for
$E^*$, equivalently if the associated twisted almost complex structure $J^*$ is linear on each fiber. We give an example showing that it is not true, which means that 
$V^*$ has no natural almost complex structure.
(Actually I believe that $J^*$ is linear only if  $J$ is integrable.)

\bigskip
Theorem 1 enables us to extend to $J$-curves in $\C\P^2$
(for a tame $J$) some classical results obtained from the theory of dual algebraic
curves. For instance, one immediately obtains the Pl\"ucker formulas, which restrict the possible
sets of singularities of $J$-curves. Such results could be interesting for the symplectic isotopy problem for surfaces in $\C\P^2$ [Sik2].
Maybe also, but this is more hypothetic, for the
topology of symplectic $4$-manifolds, 
in view of the result of D. Auroux [A] showing that they are branched coverings of $\C\P^2$. The reason why we say 
``hypothetic'' is that the branch locus is not an honest $J$-curve (with nodes and cusps), but may admit negative nodes,
although there is some some hope to dispense with them.
\bigskip

In sections 1 and 2, we give the main properties of elliptic structures and $E$-curves
in dimension $4$. In section 3 we define and study tame elliptic projective planes. Most of the statements
and all the ideas in these three sections are already in Gromov's paper (see especially [G,2.4.E and 2.4.A]), 
except what regards singularities, where we give more precise results in the vein of [McD2] and [MW].
\smallskip
In section 4 we prove the main result: the structure $E^*$ is defined in one line, but 
to prove its ellipticity we could not avoid some longish computation. The tameness and the duality property are easy.

In section 5 we give an example of a tame almost complex structure on $V=\C\P^2$ 
such that $V^*$ has no natural almost complex structure. 
\smallskip
Finally in section 6 we prove the Pl\"ucker formulas for $E$-curves and in particular for $J$-curves.
\bigskip
\noindent{\it Acknowledgment.}\hskip3mm This paper would not exist without the obstination of Stepan Orevkov, who
asked me many times  if something like an almost complex structure on $V^*$ could not exist.  I kept saying ``no'', 
until I realized the existence of an elliptic struture, which is just as good!
\bigskip
\noindent{\bf Later comment.}\hskip3mm After the first version of this paper was posted, Benjamin McKay
sent me his PhD thesis (``Duality and integrable systems of pseudoholomorphic curves'',
Duke University, 1999), which contains (among many other things) the main result of this paper in a much more general setting. 
In it one can also
find the proof of a strong version of my conjecture at the end of section 5: 
if two dual structures are almost complex near corresponding points, then they are integrable. It does not 
however contain the application to singularities.
\bigskip
\bigskip
\noindent{\bf 1. Elliptic structures in dimension $4$}
\bigskip
\noindent{\bf 1.A. Surface of elliptic type in the Grassmannian $Gr_2\R^4$.}
\hskip3mm
Let $T$ be an oriented real vector space of dimension $4$. We denote by $G(T)=Gr_2T$ the Grassmannian of oriented 
$2$-planes. Recall that for each $P\in G(T)$, the tangent 
plane $T_PG(T)$ is
canonically identified with ${\rm Hom}(P,T/P)$.
\smallskip
By definition, a {\it surface of elliptic type} in $G(T)$ is a smooth, closed, connected and embedded surface 
$E$ such that 
for every $P\in E$ one has
$$T_PE\setminus\{0\}\subset {\rm Isom}_+(P,T/P).$$
This is equivalent to the existence of (necessarily unique) complex structures on $T_PE$, $P$,
$T/P$, such that  ${\rm im}(\phi)$ is the space of complex morphisms from $P$ to $T/P$
and $\phi$ is a complex morphism. We shall denote 
\smallskip
- $j_{E,P}$ the structure on $T_PE$
\smallskip
- $j_P$ and $j_P^\perp$ the structures on $P$ and $T/P$.

\medskip
\noindent{\bf 1.B. Surfaces of elliptic type and complex structures.}
\hskip3mm The first example of surface of elliptic type is a Grassmannian $Gr_1^J(T)$ of complex $J$-lines 
for a positive complex
structure  $J$ on $T$. 

We now prove that every surface of elliptic type is deformable to such a $Gr_1^J(T)$. 
More precisely, denote by ${\cal J}(T)$ the space of positive complex structures, and ${\cal E}(T)$
the space of surfaces of elliptic type. Then the embedding ${\cal J}(T)\to {\cal E}(T)$ 
just defined admits a retraction by deformation.
 In particular, $E$ is always diffeomorphic to $\C\P^1$ and thus 
biholomorphic to $\C\P^1$.
\medskip

To prove this, we fix a Euclidean metric on $V$ and replace ${\cal J}(T)$ by the subspace ${\cal J}_0(T)$ of 
isometric structures, to which it retracts by deformation. 
The space of
$2$-vectors  $\Lambda^2T$ has a  decomposition
$\Lambda^2T=\Lambda^2_+T\oplus\Lambda^2_-T$ into self-dual and antiself-dual vectors. The 
Grassmannian $G(T)$ 
is identified
with $S^2_+\times S^2_-\subset\Lambda^2_+T\times\Lambda^2_-T$ by sending a plane $P$ to 
$(\sqrt2(x\wedge y)_+,\sqrt2(x\wedge y)_-)$ where $(x,y)$ is any positive orthonormal basis. We denote
$P=\phi(u_+,u_-)$ the plane associated to $(u_+,u_-)$. Identifying $T/P$ with $P^\perp$,
 the canonical isomorphism
$$T_{u_+}S^2_+\times T_{u_-}S^2_-\to{\rm Hom}(P,P^\perp)$$ sends $(\alpha_+,\alpha_-)$ to $A$ such that 
$$A.\xi=\iota_\xi(\alpha_++\alpha_-)$$ (interior product). This can be seen by working in a unitary
oriented
basis of $T$, $(e_1,e_2,e_3,e_4)$ such that $u_\pm={1\over\sqrt2}(e_1\wedge e_2\pm e_3\wedge e_4)$. This leads
to unitary oriented bases of $T_{u_+}S^2_+$ and $T_{u_-}S^2_-$:
$$v^\pm={1\over\sqrt2}(e_1\wedge e_3\mp e_2\wedge e_4)\,\,,\,\,w^\pm={1\over\sqrt2}(e_1\wedge e_4\pm e_2\wedge e_3).$$
Still working in these bases, one gets
 $${\rm det}\,A=-||\alpha_+||^2+||\alpha_-||^2.$$ 
(beware the signs!). Thus
an elliptic structure is given by a surface $E\subset S^2_+\times S^2_-$ such that the projections $p_\pm:E\to S^2_\pm$
satisfy 
\smallskip
- $dp_-$ is an isomorphism at all points of $E$
\smallskip
- $||dp_+\circ(dp_-)^{-1}||<1$ at all points of $E$.
\smallskip
\noindent
Since $E$ is closed and connected, $p_-$ is a diffeomorphism from $E$ to $S^2_-$ and $E$ is the set of points
$(a(u),u)$, where $a:S^2_-\to S^2_+$ is a smooth contraction.
\medskip
Thus ${\cal E}(T)$ is homeomorphic to the space  of smooth contractions from $S^2$ 
to itself. This space retracts by deformation to the space of constant maps: to each map $a$ one associates 
its unique fixed point $x_a$ and  writes 
$a=\exp_{x_h}A$ with $A:S^2\to T_{x_a}S^2$; the retraction is then $h_t=\exp_{x_h}(tA)$. Since a constant map 
$S^2_-\to S^2_+$ corresponds to a Grassmannian $Gr_1^J(T)$ for some $J\in{\cal J}_0(T)$, this 
gives the retraction by deformation from ${\cal E}(T)$ to ${\cal J}_0(T)$.
\bigskip
\noindent{\bf 1.C. Twisted complex structure associated to an elliptic surface.}\hskip3mm 
Let $E\subset G(T)$ be an elliptic
surface. Then we have the
\medskip
\noindent{\bf Proposition.} \hskip3mm {\it The space $T\setminus\{0\}$ is the disjoint union 
$\bigcup_{P\in E} P\setminus\{0\}$.}
\medskip
\noindent{\bf Proof.}\hskip3mm We use 
the representation $E=\{P_u\mid u\in S^2_-\}$, with $P_u=\phi((a(u),u)$, $a: S^2_-\to S^2_+$ being a smooth contraction.  
\smallskip
Let $\xi\in T\setminus\{0\}$ be given. Then $\xi$ belongs to $P_u$ if and only if $\iota_\xi(a(u)+u)=0$. 
We can identify  $\xi^\perp\subset T$ with $\Lambda^2_+$ in such a way that $\iota_\xi(u)=u$. 
Then $\iota_\xi(a(u)+u)=u+b(u)$, where $b:S^2_-\to S^2_-$ is a smooth contraction. Thus there exists a unique 
$u\in S^2_-$ such that
$-b(u)=u$, ie a unique $P=P_u$ in $E$ containing $\xi$. 
\smallskip
Thus if  $u$, $v$ are distinct points in $S^2_-$, we have $P_u\cap P_v=\{0\}$. In fact $P_u$ and $P_v$ are 
{\it positively}
tranvserse (first occurrence of the positivity of intersection): indeed, the inequality $||da||<1$ implies
$||a(u)-a(v)||^2<||u-v||^2$ ie 
$$\langle a(u),a(v)\rangle -\langle u,v\rangle>0,$$
which is precisely the positive transversality of $P_u$ and $P_v$.
\medskip
This enables us to put together the $j_P$, $P\in E$, to obtain a map $J:T\to T$ with the following properties:
\smallskip
\item{} (i)$J^2=-{\rm Id}$
\smallskip
\item{} (ii) $J$ is continuous, and homogeneous of degree $1$
\smallskip
\item{} (iii) $J$ is smooth away from $0$ (by (ii), it is not differentiable at $0$ except if it is linear)
\smallskip
\item{} (iv) for every $x\in T\setminus\{0\}$, $J(x)$ is linearly independent of $x$, and $J$ is linear on the plane 
$\langle x,J(x)\rangle$.
\smallskip
\noindent Conversely, given $J$ satisfying (i)-(iv), we can define a smooth surface $E\subset G(T)$ by 
$$E=\{\langle x,J(x)\rangle\mid x\in T\setminus\{0\}\}.$$ A straightforward computation gives that
$E$ is elliptic if and only if
\smallskip
\item{} (v)
for every $x\in T\setminus\{0\}$  and $\xi\in \langle x,J(x)\rangle^\perp\setminus\{0\})$,
$(x,J(x),\xi,dJ_x.\xi)$ is an oriented basis of $T$.
\bigskip
\noindent{\bf 1.D. Local form of a surface of  elliptic type.}
\hskip3mm Let $E$ be a surface of elliptic type in $G(T)$, and fix $P\in E$. Identify $T$ with $\C^2$ 
such that
\smallskip
(i) $P$ is sent to the horizontal plane $H=\C\times\{0\}$
\smallskip
(ii) the identifications $P\to H$ and $T/P\to\C^2/H$ are complex-linear for the complex structures defined in 1.A.
\smallskip
\noindent Then $E$ is given near $P$  by a 
family of
planes of the following form:
$$P_{\lambda}=
\{(\delta z,\delta w)\in\C^2\mid \delta w=\lambda\delta z+h(\lambda)\overline{\delta z}\},$$
where $h$ is a smooth germ $(\C,0)\to(\C,0)$ such that $h(0)=0$. The tangent space $T_{P_\lambda}E$ is identified
with the image of
$$\xi\in\C\mapsto\xi{\rm Id}+(dh(\lambda).\xi)\sigma\in{\rm End}_{\R}(\C),$$ 
where $\sigma$ is the complex conjugation. Thus the ellipticity translates to the inequality $||dh||<1$. Property
(ii) becomes $dh(0)=0$.
\medskip
\noindent{\it Remark.} Denote by $J$ the complex structure on $T$. Properties (i) and (ii) imply that 
$P\in Gr_1^JT$ and $T_PGr_1^J(TV)=T_PE$, cf. [G, 2.4.E].

\bigskip
\noindent{\bf 1.E. Elliptic structure on a $4$-manifold; twisted almost complex structure.}\hskip3mm
Let $V$ be an oriented $4$-manifold, and let $G=Gr_2TV$ be the Grassmannian of oriented tangent $2$-planes, which is
fibered over $V$ with fiber $G_v=Gr_2(T_vV)$.
\smallskip
By definition, an {\it elliptic structure} on  $V$ is a smooth compact
submanifold $E\subset G$ of dimension $6$, transversal to the fibration $G\to V$, such that
each fiber  $E_v$ is a surface of elliptic type in $Gr_2T_vV$.
\smallskip
Denote ${\cal E}(TV)$ the space of  elliptic structures on $V$ and ${\cal J}(TV)$
the space  of positive almost complex structures on $V$. The map $J\mapsto Gr_1^J(TV)$ gives a
 natural embedding from ${\cal J}(TV)$ to ${\cal E}(TV)$. These are both spaces of sections of 
a bundle on $V$, with respective fibers
${\cal E}(T_vV)$ and ${\cal J}(T_vV)$. Since ${\cal E}(T_vV)$ retracts by deformation to ${\cal J}(T_vV)$, 
${\cal E}(TV)$ retracts by deformation to ${\cal J}(TV)$. In particular, every elliptic structure defines a
unique homotopy class of almost complex structures on $V$. Thus the Chern class 
$c_1(E)=c_1(TV,J)\in H^2(V,\Z)$ is well defined.
\medskip
Finally, the twisted structures $J_v$, $v\in V$, can be put together to give a 
{\it twisted almost complex structure}
on $V$, ie a fiber-preserving map $J:TV\to TV$ such that all the $J_v$ have the properties (i)-(v) of 1.C. 
It is continuous on $TV$ [in fact locally Lipschitz], and smooth away from the zero section. Conversely, a map $J$ with all 
these properties clearly
defines an elliptic structure.

\bigskip
\bigskip

\noindent{\bf 2. Solutions of $E$, $E$-maps and $E$-curves}
\bigskip
In this section we consider an oriented $4$-manifold $V$ equipped with an elliptic structure $E\subset G=Gr_2TV$.
 If $S$ is an oriented surface and 
$f:S\to V$ is an immersion, we denote $\gamma_f:S\to G$ the associated Gauss map.
\bigskip
\noindent{\bf 2.A. Immersed solutions. Local equation as a graph }

\medskip
\noindent{\bf Definition.} An immersed {\it solution of $E$} is a ${\cal C}^1$ immersion $f:S\to V$ where $S$ 
is an oriented surface and $\gamma_f(S)\subset E$.
\smallskip
Let $v\in V$ and $P\in E_v$ be fixed. We describe a local equation for germs of
immersed solutions of $E$ which are tangent to $P$ at $v$, 
or more generally which have a tangent close enough to $P$.
Choose a local chart $(V,v)\to(\C^2,0)$ such that the properties of 1.D are satisfied for $E_v$ and $P$. 
Then the elliptic structure on $E$ near $H$ is given by a 
family of
planes
$$P_{z,w,\lambda}=
\{(\delta z,\delta w)\in\C^2\mid \delta w=\lambda\delta z+h(z,w,\lambda)\overline{\delta z}\}.$$
Here $(z,w,\lambda)$ belongs to a neighbourhood of $0$ in $\C^3$, and
$P_{z,w,\lambda}$ represents a plane tangent at the point of coordinates $(z,w)$. The map
$h$ is a smooth germ $(\C^3,0)\to(\C,0)$ such that
$$\left\{\eqalign{
&||D_3h(z,w,\xi)||<1\hskip3mm (\forall (z,w,\xi))\cr
&D_3h(0,0,0)=0.\cr}\right.\leqno(1)$$
A germ of surface $S\subset V$ passing through $P$ with a tangent plane close enough to $P$, can be written
as a graph $w=f(z)$, where $f:(\C,0)\to(\C,0)$ satisfies 
$${\partial f\over\partial\overline z}=h(z,f(z),{\partial f\over\partial z}).\leqno(2)$$
\noindent{\it Remark.} This equation  with the property $||D_3h||<1$ is the general (resolved) form of 
an elliptic equation $\C\to\C$, cf. [V]. It implies the existence of a local immersed solution of $E$ with 
any given tangent plane, and even with an arbitrary ``compatible'' $k$-jet (an easy proof can be given by a suitable 
implicit function theorem, modifying slightly
the proofs given in Chapter V or VI of [AL]), and also that each solution is of class
${\cal C}^\infty$.

\bigskip
\noindent{\bf 2.B. Conformal parametrization, $E$-maps.}\hskip3mm Let $f:S\to V$ be an immersed solution of $E$.
Since every tangent plane $P_z=df_z(T_zS)$ has a well-defined complex structure 
$j_{P_z},$ this induces a canonical 
almost complex structure
$j_f$ on $S$, ie a natural structure of Riemann surface. In other words, every (immersed) solution of $E$ admits
a natural conformal parametrization. 
\smallskip
If $S$ is a Riemann surface, we say that an immersion $f:S\to V$ is a 
{\it conformal solution} of $E$
if it is a solution and $j_f$ is the canonical almost structure on $S$. This is equivalent to the equation
$$df_z\circ i=J_{f(z)}\circ df_z.\leqno(3)$$
We can now eliminate the immersion condition and define an {\it $E$-map} as a ${\cal C}^1$ map $f:S\to V$, where 
$S$ is a Riemann surface, which is
a solution of (3). 
\smallskip
Note that since $J$ is only Lipschitz, the fact that $E$-maps are smooth is not completely obvious at this stage. 
But the arguments of [AL, chap. V or VI] imply that if $f$ is a nonconstant local $E$-map, then
$df$ has only isolated zeros and the Gauss map $\gamma_f$ can be extended continuously at these zeros.
And also that there exist $E$-immersions with an arbitrary given $2$-jet, 
satisfying suitable compatibility conditions.
\bigskip
\noindent{\bf 2.C. $E$-maps as pseudoholomorphic maps.}\hskip3mm We prove here that $E$-maps 
can be considered as pseudoholomorphic maps. For this, we define
the $4$-dimensional holonomy distribution $\Theta$ on $E$ by setting  $\Theta_P= d\pi_P^{-1}(P)$, where 
$\pi:G\to V$ is the natural
projection. It is characterized by the property that 
every
Gauss map $\gamma_f:S\to E$ associated to an $E$-map (not locally constant), is 
tangent to $\Theta$.  The construction of [G],1.4 can be generalized to give the
\medskip
\noindent{\bf Proposition.}\hskip3mm {\it There exists a unique almost complex structure $\widetilde J$ on $\Theta$, 
such that
every Gauss map $\gamma:S\to E$ associated to an $E$-map is $\widetilde J$-holomorphic 
(or is a $\widetilde J$-map), ie it is ${\cal C}^1$ (or ${\cal C}^\infty$),
tangent to $\Theta$ and satisfies
$$d\gamma_z\circ i=\widetilde J_{\gamma(z)}\circ d\gamma_z.\leqno(4)$$
Moreover,  the differential $d\pi_P:T_PG\to P$ is complex linear on $\Theta_P$
and the subbundle $F$ tangent to the fibers is $\widetilde J$-invariant, and ${\widetilde J}_{P|F_P}=j_{E_v,P}$ (notation of 1.A).
\smallskip
Conversely, a $\widetilde J$-map not locally contained in a fiber is the Gauss map associated to an $E$-map. 
}
\medskip
\noindent{\bf Sketch of proof.}\hskip3mm For every
$X\in \Theta_P\setminus F_P$ there is an $E$-map $f:S\to V$ and a vector $u\in T_zS$ such that $d\gamma_z(u)=X$
where $\gamma$ is the Gauss map. Thus necessarily
$\widetilde J_P.X=d\gamma_z(iu),$
which implies the uniqueness.

\smallskip
The existence of $\widetilde J$ can be proved in local coordinates. 
The equation (3) becomes
$${\partial f\over\partial y}=J_{f(z)}({\partial f\over\partial x}),\leqno(5)$$ 
and the associated Gauss map is $\gamma=(f,s)$ with
$$s(z)={{\partial f\over\partial x}\wedge(J_{f(z)}({\partial f\over\partial x}))\over
||{\partial f\over\partial x}\wedge(J_{f(z)}({\partial f\over\partial x}))||}.$$
By differentiating (3) and after some (straigthforward but tedious) computation, one can obtain 
an equation of the form 
${\partial\gamma\over\partial y}=\widetilde J.{\partial\gamma\over\partial x}$ and check that $\widetilde J$
has all the stated properties. The uniqueness implies that it is independent of the chart.

\smallskip
We know already that $\gamma$ is ${\cal C}^\infty$ away from singularities of $f$. If $z$ is such 
a singularity, then since $\gamma$ satisfies (4) away from $z$ and is continuous up to $z$, 
it is ${\cal C}^\infty$ up to $z$. Note that it implies that every $E$-map is ${\cal C}^\infty$.
\smallskip
Finally, if $\gamma:S\to E$ is a ${\cal C}^1$  map, tangent to $\Theta$ and satisfying (4),
 the composed map $f=\pi\circ\gamma$ satisfies (3). If $\gamma$ is not locally contained in a fiber,
$f$ is not locally constant, thus its Gauss map $\gamma_f$
is well-defined and one has $\gamma_f=\gamma$.
\medskip
\noindent{\bf Remarks}
\smallskip (i) One can view $\widetilde J$-map as an ordinary 
pseudoholomorphic 
map by extending
$\widetilde J$ arbitrarily to $\widehat J$ defined on $TE$ (using  Riemannian metric for instance). Then
a $\widetilde J$-map is a $\widehat J$-map which is tangent to $\Theta$.

\smallskip
(ii) If $F$ was totally real instead of being complex, the triple $(\Theta,J,F)$ 
would be
(adding some partial integrability
properties)  what F. Labourie [L] calls 
a {\it Monge-Amp\`ere geometry}. The main difference is that his compactness theorem requires
to add ``curtain surfaces'', ie $\widehat J$-curves satisfying ${\rm dim}(TS\cap F)\equiv 1$.

\bigskip
\noindent{\bf 2.D. General $E$-curves, compactness theorem.}\hskip3mm Using the conformal parametrization, we can now define an 
$E$-curve
as a $\widetilde J$-curve $\widetilde C$ in $E$
which is ``almost transverse'' to $F$. To make this definition precise, we would have 
\smallskip
(i) to choose a definition of $\widetilde J$-curve, 
for instance a stable $\widetilde J$-curve in the sense of Kontsevich. 
\smallskip
(ii) to say what ``almost transverse'' means: essentially that no nonconstant
component has an image contained in a fiber, or that it is transverse to $F$ except on a finite subset.
\medskip
We would then obtain topological spaces of $E$-curves. We shall not give any details here since the only
spaces of $E$-curves we shall consider will consist of curves which are embedded (in $V$!). The
 space of embedded $E$-curves will be considered as a subspace of the space of smooth surfaces in $V$: 
recall that this is a Fr\'echet
smooth manifold whose tangent space at $S$ is the space of sections of the normal bundle $N(S,V)=T_SV/TS$.
\medskip
We shall also use the following notions of individual $E$-curves:
\smallskip
- a {\it primitive} (or irreducible) $E$-curve is the image $C=f(S)$ where $S$ is a closed and connected Riemann
surface, $f$ is an $E$-map which does not factor $f=f_1\circ\pi$ with $\pi$ a nontrivial holomorphic covering. 
As in the case of $J$-curves, the image determines $(S,f)$ up to isomorphism
\smallskip
- an {\it $E$-cycle} (\`a la Barlet) $C=\sum n_iC_i$ where the $C_i$ are distinct
primitive $E$-curves and the $n_i$ are positive integers
\smallskip
- analogous local versions of these.
\bigskip
One expects a compactness theorem for $E$-curves, analogous to the one for pseudoholomorphic curves: 
roughly speaking, it should say that (if $V$ is compact) a set of  $E$-curves is relatively compact
if their areas in $V$ are uniformly bounded. If one replaces ``areas in $V$'' by ``areas in $E$ of the Gauss
maps'', then such a result follows from 
\smallskip
- the  compactness theorem for $\widehat J$ on $TE$ (cf. the remark (i) in 2.C)

- the fact that the conditions ``tangent to $\Theta$'' and ``almost transverse to $F$'' are closed conditions.
\smallskip

However, it is not clear  that an area bound in $V$ gives an area bound in $E$. Gromov [G, 2.4.E] says that
the Schwarz lemma is still valid for $E$-curves under an area bound in $V$, but I do not understand  the proof.
Anyhow, here  we shall only need the compactness theorem for $E$-lines, which we shall prove in section 3.
\bigskip
\noindent{\bf 2.E. Singularities  of $E$-curves and positivity of intersections}
\medskip
Here we extend to $E$-curves the result of M. Micallef and B. White
[MW] (see also [Sik1]): we prove that a $E$-curve, possibly non reduced, is ${\cal C}^1$-equivalent to a germ of standard
holomorphic curve in $\C^2$. It implies the positivity of intersections for $E$-curves, 
in particular
the genus and intersection formulas.
\medskip
Such a result could be proved by showing that such a surface is {\it quasiminimizing} in the sense 
of [MW], but we prefer to use more complex-analytic arguments as in [Sik1]. 
\medskip
We use the chart of 2.A to write the equation in intrinsic form, ie as a graph over the tangent space. 
If the curve is non 
singular, this is the equation (2) where $h$ satisfies (1). Now, consider a germ of non-immersed $E$-map 
$F:(\C,0)\to(\C^2,0)$ with horizontal tangent at the origin. Then  equation (3) and the similarity principle [Sik1]
(proposition 2; cf. also [McD2])
give
the existence of $a\in\C^*$ and $k\in\N$, $k\ge2$, such that
$$f(z)=(az^k,0)+O_{2,1-}(z^{k+1}).$$
Here we  use some notation from [MW], [Sik1]: $g(z)=O_{2,1-}(z^{k+1})$ means

\smallskip
$$\left\{\matrix{g(z)=O(z^{k+1}),dg(z)=O(z^k),d^2g(z)=O(z^{k-1})\cr
(\forall\alpha<1)\,d^2g \,\hbox{\rm is}\,\alpha\hbox{\rm -H\"older  with 
H\"older constant}\,\, O(|z|^{k-1-\alpha}).\cr}\right.$$
\smallskip
\noindent 
Thus we can reparametrize the curve by setting $pr_1\circ F(z)=t^k$, where $z\mapsto t$ 
is a ${\cal C}^1$ local diffeomorphism. We obtain a map $t\mapsto (t^k,F(t))$ where $F$ is of 
class ${\cal C}^{2,1-}=\cap_{\alpha<1}{\cal C}^{1,\alpha}$, with $F(t)=O_{2,1-}(t^{k+1})$
(cf. [Sik1], proposition 4). The fact 
that the image, viewed locally as a graph satisfies (1), means that 
$F$ satisfies
$${\partial F\over\partial\overline t}=k\overline t^{k-1}\,h(t^k,F(t),{1\over kt^{k-1}}
{\partial F\over\partial t}).\leqno(6)$$
Note that this makes sense also at the origin. 
\medskip
Next, we show that
the ``difference'' of two such germs satisfies the similarity principle.
More precisely,  we have the following generalization of [Sik1], prop. 5:
\medskip
\noindent{\bf Proposition.}\hskip3mm {\it Assume that $F$ and $G$ satisfy (6) 
with the same value of $k$, and
are not identical as germs. Then there exists $a\in\C^*$
and $k\in\N^*$ such that $$F(t)-G(t)=at^k+O_{1,1-}(t^{k+1}).$$ }
\noindent{\bf Proof.}\hskip3mm  Set $u=F-G$, and take the difference of the two equations 
on $F$ and $G$. Using Taylor's integral formula and setting
$$\gamma(t,s)=\big(t^k,G(t)
+su(t),{1\over kt^{k-1}}({\partial G\over\partial t}+s{\partial u\over\partial t})\big),$$ we get
$$\eqalign{{\partial u\over\partial\overline t}&=A(t).u(t)+B(t).{\partial u\over\partial t},\cr}$$
where
$$\left\{\eqalign{A(t)&=\int_0^1k\overline t^{k-1}D_2h(\gamma(t,s))\,ds\cr
B(t)&=
\int_0^1 \overline kt^{k-1}D_3h(\gamma(t,s))\,ds.\cr}\right.$$
The properties of $h$, $F$ and $G$ imply that $A$ and $B$ are of class 
${\cal C}^{1,1-}$, and 
$||B||_{L^\infty}<1$. Then the proposition follows from a variant of proposition 2 in [Sik1].
\medskip
Finally, one proceeds exactly as in [Sik1] (inspired by [MW]) to deduce from this proposition the 
\medskip
\noindent{\bf Proposition.}\hskip3mm {\it Let $E$ be a germ of elliptic structure on $\C^2$ 
near $0$
such that the horizontal plane
$H=\C\times\{0\}$ belongs to $E_0$.
Let $f_i:(\C,0)\to (\C^2,0)$, $i=1,\cdots,r$, be germs
of $E$-maps, all tangent to $H$ at $0$.
Then there exist 

\smallskip 
- a local ${\cal C}^1$ diffeomorphism $\phi:(\C^2,0)\mapsto(\C^2,0)$, with support in an arbitrarily
small sector 
$$S_\epsilon=\{(x,y)\in\C^2\mid|y|\le\epsilon|x|\}$$

- local diffeomorphisms $u_i:(\C,0)\mapsto(\C,0)$, tangent to the identity,
\smallskip
\noindent 
such that all the maps 
$\phi\circ f_i\circ u_i$ are holomorphic.
} 
\medskip
If the tangents to the $f_i$ are not the same, we cannot in general expect to find a differentiable
chart on $V$ in which the image become 
holomorphic: there is an obstruction already at the linear algebraic level. However, by superposing the diffeomorphisms given by the proposition
 we easily obtain a Lipschitz chart:
\medskip
\noindent{\bf Theorem.}\hskip3mm {\it Let $f_i:(\C,0)\to (V,v)$, $i=1,\cdots,r$, be germs
$E$-maps through the same point.  Then there exists a germ of Lipschitz oriented homeomorphism 
$\phi:(V,v)\to(\C^2,0)$
 such that all the maps 
$\phi\circ f_i\circ u_i$ are holomorphic.}
\medskip
\noindent{\bf Proof.}\hskip3mm We can assume that $(V,v)=(\C^2,0)$. Let
$(P_j)$, $j=1,\cdots,r$, be the different tangent planes to the $f_i$ at $v$, and let
$I_j\subset\{1,\cdots,r\}$ be the indices corresponding to the branches
with tangent $P_j$. The $P_j$ are not
complex linear in general, but there exist a Lipschitz oriented homeomorphism $h$ of $\C^2$ such that the $h(P_j)$
are complex linear, thus there exist a family of complex linear $A_j$ such that
 $$(\forall j)\,\,A_jh(P_j)=\C\times\{0\}.$$ Furthermore, 
we can assume that $h$ is smooth [even linear] on a sector $S_j$ around $P_j$.

\smallskip
Then we can apply the proposition to $A_jhf_i$, $i\in I_j$: there exists a ${\cal C}^1$ diffeomorphism
 $\phi_j:\C^2\to\C^2$, such that $\phi_jA_jhf_i$ is holomorphic
with horizontal tangent for $i\in I_j$. Moreover, we can assume that the support of $\phi_j$
is contained in $A_jh(S_j)$.
\smallskip
Then the desired homeomorphism is given by
$$\phi=\left\{\matrix{A_j^{-1}\circ \phi_j\circ A_j\circ h\,\,{\rm on}\,\,S_j\cr
h\,\,{\rm elsewhere}.\cr}\right.$$
\medskip
From this theorem, one deduces the positivity of intersections. More precisely, one can define
\smallskip
- a local intersection index $(C,C')_v\in\N^*$ for two germs
$C$ and $C'$ of $E$-cycles at the point $v$ without common component. It is equal
to $1$ if and only if $C$ and $C'$ are smooth at $v$, with distinct tangents
\smallskip
- a local self-intersection number $\delta_v(C)\in\N^*$ (number of double points in a generic
deformation equitopological at the source). It is equal
to $1$ if and only if $C$ is smooth at $v$.
\medskip
One then has the intersection and genus (or ``adjunction'') formulas
\medskip
\noindent{\bf Theorem.}
\smallskip {\it (i) If $C$ and $C'$ are two $E$-cycles without common components, then $C\cap C'$ 
(intersection of the supports)
is finite and
the homological intersection is given by
$$C.C'=\sum_{v\in C\cap C'}\,(C,C')_v.$$ 

(ii) If $C$ is an irreducible $E$-curve of genus $g$, then it has a finite number of singularities, and 
its genus is given by
$$g={C.C-c_1(E).C\over2}+1-\sum_v\,\delta_v.$$} 

Assume that $(V,E)$ is homologically equivalent to $(\C\P^2,E_0)$. Thus there is a well
defined degree map $H_2(V;\Z)\to\Z$ (an isomorphism modulo torsion, 
but not necessarily an isomorphism at this stage). Define an {\it  $E$-line} as a primitive $E$-curve $C\subset V$ of degree $1$. Then 
Since an $E$-line satisfies $C.C=1$ and $c_1(E).C=3$, we get the
\medskip
\noindent{\bf Corollary.}\hskip3mm{\it Every $E$-line is  an embedded sphere, and two distinct 
$E$-lines intersect transversely in one unique point.}
\bigskip

\noindent{\bf 2.F. Linearization of the equation of $E$-curves; automatic genericity.}\hskip3mm
We consider here embedded $E$-curves, ie smooth surfaces $S\subset V$ satisfying $\widetilde S\subset E$
where $\widetilde S$ is the Gaussian lift. Following [G, 2.4.E]
we linearize this ``equation'' at $S$, obtaining an equation $\overline\partial_E f=0$ where $\overline\partial_E$ 
acts on sections $f:S\to N=T_SV/TS$, the normal bundle, with values in $\Omega^{0,1}_J(S,N)$. 
Here $J=j^\perp_{TS}$ (cf. 1.A) is the natural complex structure on $N$.
\smallskip
One can obtain explicitly this equation by using the equation (2) in local coordinates. The compatibility of the complex structure
on $\C^2$ with the structure on $N$ means that $D_3h(z,0,0)\equiv0$. Thus the linearization of (2) has locally the form
$${\partial f\over\partial\overline z}-D_2h.f=0,$$
ie the operator $\overline\partial_E$ has the form $\overline\partial_E=\overline\partial+R$ where 
$\overline\partial$ is associated to a holomorphic structure on $N$ and $R$ is of order $0$.
Thus one can apply to it the arguments 
of [G, 2.1.C] (cf. also [HLS]): 
\medskip
\noindent{\bf Proposition.} \hskip3mm{\it If $c_1(N)>2g-2$, ie $c_1(E).S>0$, then $\overline\partial_E$ is onto. 
Thus the space $M_A$ of connected embedded $E$-curves in the class $A\in H_2(V;\Z)$, if nonempty,
 is a smooth manifold if $c_1(E).A>0$. 
Its real dimension is 
$$2(c_1(N).S+1-g)=2A.A+(c_1(E).A-A.A)=A.A+c_1(E).A.$$}

Also, $M$ is oriented since the homotopy $\ker (\overline\partial+R_t)$ gives it a natural
homotopy class of almost complex structure (not more than a homotopy class, as we prove in section 5!).
\medskip
Now assume that $(V,E)$ is homologically equivalent to $(\C\P^2,E_0)$. Then the space of $E$-lines is the
disjoint union
of the $M_A$ for all $A\in H_2(V;\Z)$ of degree $1$. 
For such classes, we have $c_1(E).A=3$ thus the proposition applies:
\medskip
\noindent{\bf Corollary.}\hskip3mm {\it The space of 
$E$-lines $V^*$, if nonempty,  is naturally a smooth oriented $4$-manifold.} 
\medskip 
\noindent{\bf Extensions.}\hskip3mm Since $c_1(E).A=3$ we can still impose on $S$ a condition of complex codimension $1$ or $2$, and
keep the automatic genericity (cf. for instance [B]). In particular:
\smallskip
(i) Let ${\cal L}^*_v$ be the space of $E$-lines through a given $v\in V$: it is an oriented smooth surface in  $V^*$ when nonempty. Note that
it can be identified
with an open subset of the projective line $G_1^J(T_vV)$. Also,
${\cal L}^*_v$ depends
smoothly on $v$.
\smallskip
(ii) Let ${\cal L}^*_{v,w}$ be the space of $E$-lines through two given points $v, w\in V$: 
when nonempty, it is a point $L_{v,w}$ which depends smoothly on $(v,w)$. This is the case for some open subset 
$U_1\subset V\times V\setminus\Delta_V$.
\smallskip
(iii) Let ${\cal L}^*_P=$ be the space of $E$-lines with a given tangent plane $P\in E$: 
again, when nonempty, it is a point $L_P$ which depends smoothly on $P$. This is the case for some open subset
$U_2\subset E$.
\bigskip
\bigskip
\noindent{\bf 3. Tame elliptic projective planes}
\medskip
By definition, a {\it tame elliptic projective plane} $(V,E)$ is
a $4$-manifold equipped with a tame elliptic structure,  
homologically equivalent to $(\C\P^2,E_0)$. Note that we do not require a priori $V$ to be diffeomorphic to $\C\P^2$.
\bigskip
\noindent{\bf 3.A. Proposition.}\hskip3mm{\it Let $(V,E)$ be a tame elliptic projective plane. Then 
\smallskip
(i) by two distinct
points $x,y\in V$ there passes a unique $E$-line $L_{v,w}$, and for any given 
$P\in E$ there exists a unique $E$-line $L_P$ tangent to $P$.
\smallskip
(ii) $V$ is oriented diffeomorphic to $\C\P^2$, $E$ is homotopic to $E_0$, and
any taming $\omega$ is isomorphic to $\lambda\omega_0$ for some $\lambda>0$.}
\medskip
\medskip
\noindent{\bf Proof} 
\smallskip 
(i) It suffices to prove that  $V^*$ is compact: this will imply that that the open sets  
$U_1$ and $U_2$ defined at the end of 2.F
are also closed, so 
$U_1= V\times V\setminus\Delta_V$ and $U_2=E$, which proves (i).

\smallskip
The compactness of $V^*$ will follow
from the compactness theorem for $\widehat J$. First, there exists a taming $\Omega$ for $\widehat J$: 
as usual in the theory of symplectic bundles, we set
$\Omega=\pi^*\omega+\alpha$ where $\alpha$ is a closed $2$-form on $E$ which is positive on very fiber. 
Such a form exists since $H^2(E)\to H^2(V)$ is onto: this is true since it holds for in the standard case
$V=\C\P^2$, $E=E_0$ and that our case
is homologically standard. 
\smallskip
Furthermore, let $A\in H_2(E;\Z)$ be the homology class of the
Gauss lift of $E$-lines. Then $A$ is {\it $\Omega$-indecomposable}, ie not equal to a sum $A=A_1+A_2$
with $\omega(A_i)>0$. This can  again be seen in the standard situation 
[in that case, a holomorphic curve $C$ in the class $A$ is always a section $s(L)$ over a line $L$ in $\C\P^2$; 
if $C$ is not the Gaussian lift of $L$, then there exists $v_0\in \C\P^2\setminus L$ such that $s(v)$ is the tangent to
the
line
$[v_0v]$ for every $v\in L$].
\smallskip
The $\Omega$-indecomposability and the compactness theorem of [G] imply that the space $M$ of  rational $\widehat J$-curves in the class $A$ is compact, and since $V^*$
is homeomorphic to a closed subset of $M$ it is also compact.
\medskip
(ii) Fix three lines $L_0$, $L_1$, $L_\infty$. We deform $E$ so that it remains tamed by $\omega$, the $L_i$ are still $E$-lines and 
$E$ comes from
 a complex 
structure isomorphic to the standard one near $L_\infty^0$. This is possible, using Darboux-Givental and
the contraction ${\cal E}_\omega\to{\cal J}_\omega$.

\smallskip
We shall find a diffeomorphism $\phi:V\to\C\P^2$ which sends
them to the $x$-axis $L_0^0$, the $y$-axis $L_1^0$ and the line at infinity $L_\infty^0$.

\smallskip Denote  $v_0$, $v_1$ the intersections 
$L_0\cap L_\infty$, $L_1\cap L_\infty$. 
Let $v\in V\setminus L_\infty$. Then the $E$-lines $v_0v$ and $v_1v$ meet $L_0$ and $L_1$ in 
$x(v)$ and $y(v)$ respectively. Identifying $L_0$ with $L_0^0$, $L_1$ with $L_1^0$, we define $\phi(v)$ to be the 
intersection of $v_0x(v)$ and $v_1y(v)$. We obtain thus a smooth map
$\phi:V\setminus L_\infty\to \C\P^2\setminus L_\infty^0$
\smallskip
Exchanging the roles of $V$ and $\C\P^2$, we obtain $\psi:V\setminus L_\infty^0\to \C\P^2\setminus L_\infty$ which 
is the inverse of $\psi$.
Since everything is standard near $L_\infty$, one can extend $\phi$ to $L_\infty$ and $\psi$ to $L_\infty^0$.
\medskip
Finally, the fact that $\omega$ is isomorphic to $\lambda\omega_0$ results from Moser's lemma.
\bigskip
\noindent{\bf 3.B.}\hskip3mm Conversely, as shown by Gromov, one has the
\medskip
\noindent{\bf Proposition  [G, 2.4.A'].} \hskip3mm {\it Assume that $V^*$ is compact and nonempty. 
Then there exists a taming symplectic form $\omega$.}

\medskip 
\noindent{\bf Proof.}\hskip3mm Using a positive volume form $\nu$
on $V^*$ (identified with
a smooth measure $d\nu$), define a $2$-form $\omega$ by Crofton's formula:
$$\int_S\omega=
\int_{V^*}\,{\rm Int}\,(S,L)\,d\nu(L)$$
for every oriented surface $S\subset V$. Here, ${\rm Int}(S,L)$ is the algebraic intersection number, which is defined for
almost all $L\in V^*$.
\smallskip
Let us give a more explicit definition of $\omega$. First, fix $v\in V$ and denote by ${\cal L}^*_v\subset V^*$ 
the subset of $E$-lines containing a given $v\in V$, which is a submanifold 
diffeomorphic to 
$\C\P^1$.

\medskip

\noindent{\bf Proposition.} \hskip3mm {\it There is a canonical isomorphism between normal bundles, $\nu:N_vL\to 
N_L{\cal L}^*_v$.}
\medskip
\noindent{\bf Proof.} 
\hskip3mm Choose an $E$-line $L^\perp$ different from $L$ at $v$. If 
$\delta L\in T_LV^*$, let $(L_t)$ be a path such that ${d\over dt}_{t=0}L_t=\delta L$, and set
$$\phi(\delta L)=
\big({d\over dt}\big)_{|t=0}(L_t\cap L^\perp)\,\in T_vV.$$
Dually, choose a point $w\in L$ different from $v$. If 
$\delta v\in T_vV$, let $(v_t)$ be a path such that ${d\over dt}_{t=0}v_t=\delta v$, and set
$$\psi(\delta v)=
\big({d\over dt}\big)_{|t=0}L_{v_t,w}\,\in T_LV^*.$$
Then clearly, $\phi$ induces the desired isomorphism $\nu$ and $\psi$ its inverse.
\medskip
We can now define a morphism
$$s:T_vV\to\Gamma({\cal L}^*_v,N({\cal L}^*_v,V^*))$$
 by composing
$$T_vV\to N_vL\to N_L{\cal L}^*_v\,\,,\,\,L\in {\cal L}^*_v.$$ Thus for
$X, Y\in T_vV$ and $L\in V^*_x$, $s(X)(L)$ and $s(Y)(L)$ are elements
of $N_L{\cal L}^*_v$. Lifting them
to $\widetilde X,\widetilde Y\in T_LV^*$, we see that $$\nu(\widetilde X,\widetilde Y)=\iota_{\widetilde X}\iota_{\widetilde Y}\nu\in
\Lambda^2T^*_LV^*$$
(interior products) is independent of the lifts. Varying $L\in {\cal L}^*_v$, we obtain a $2$-form on ${\cal L}^*_v$
which we denote by $\nu(s(X),s(Y))$, and we set

$$\omega(X,Y)=\int_{{\cal L}^*_v}\nu(s(X),s(Y)).$$
It is easy to see that it is positive on $E$ and satisfies Crofton's formula.
\bigskip
\bigskip
\noindent{\bf 4. Proof of the main result}
\medskip
\noindent{\bf 4.A. Definition of $E^*$.}\hskip3mm We set
$$E^*=\{T_L{\cal L}^*_v\mid v\in V, L\in {\cal L}^*_v\}.$$
This is clearly a submanifold fibered over $V^*$, the fiber at $L$ being
$$E^*_L=\{T_L{\cal L}^*_v\mid v\in L\}.$$
It is equipped with a natural distribution of codimension $2$, $\Theta_{P^*}=d\pi^{*-1}(P^*)$.
Note also that $E^*$ is naturally diffeomorphic to $E$ via $\phi:T_vL\mapsto T_L{\cal L}^*_v$, 
in fact both are naturally diffeomorphic to the 
incidence variety
$$I=\{(v,L)\in V\times V^*\mid v\in L\}=\{(v,L)\in V\times V^*\mid L\in {\cal L}^*_v\}.$$
This variety is equipped with two natural fibrations $p:I\to V$, $p^*:I\to V^*$. 
It also has one natural distribution. Indeed, by differentiating the condition $(v(t)\in L(t))$, one obtains the
\medskip
\noindent{\bf Proposition.} \hskip3mm {\it If $(\delta v,\delta L)\in T_{v,L}I$, then
$(\delta v\in T_vL\Leftrightarrow\delta L\in T_L{\cal L}^*_v).$}
\medskip
Thus one can define the distribution $D\subset TI$ by
$$D_{v,L}=dp^{-1}(T_vL)=dp^{*-1}(T_LL_v^*).$$
We then have a commutative triangle
$$\matrix{(I,D)\cr
\noalign{\vskip 1mm}
\gamma\swarrow\hskip10mm\searrow\gamma^*\cr
\noalign{\vskip 3mm}
\diagram
(E,\Theta)&\hfl[0mm][9mm][0mm]{\phi}{}&(E^*,\Theta^*)\cr
\enddiagram}$$
where $\gamma(v,L)=T_vL$ and $\gamma^*(v,L)=T_LL_v^*$.
\bigskip
\noindent{\bf 4.B. Proof that $E^*$ is elliptic}\medskip 
Let us fix $(v,L)$ such that $v\in L$, and denote 
$$P=T_vL\,\,,\,\,P^*=T_L{\cal L}^*_v.$$
Then there are natural embeddings
$$\eqalign{i&:T_PE_v\to{\rm Hom}(T_vL,N_vL)\cr
i^*&:T_{P^*}E^*_L\to{\rm Hom}(T_L{\cal L}^*_v,N_L{\cal L}^*_v).\cr}$$
By the ellipticity of $E$, $i(p)$ is an oriented isomorphism if $p\ne 0$. We want to prove the same property for
$i^*(p^*)$. This will follow from
\smallskip
- the existence of canonical isomorphisms $T_PE_v\approx P^*$, $:T_{P^*}E^*_L\approx P$, $N_vL\approx N_L{\cal L}^*_v$
(this last we know already); thus $i$ and $i^*$
become morphisms $P^*\to{\rm Hom}(P,N)$ and $P\to{\rm Hom}(P^*,N)$
\smallskip
- the formula
$$i^*(p)(p^*)=i(p^*)(p).$$
To prove this, we define local charts on $V$, $V^*$ and $G$: 
\smallskip
1) We start with a chart $\Psi:V\to T_vL\times N_vL$, such that
$$\left\{\eqalign{&\Psi(L)=T_vL\times\{0\}\cr
&\Psi(L^\perp)=\{0\}\times N_vL\cr
&{\rm pr}_1\circ d\Psi_{v|T_vL}={\rm Id}\cr
&{\rm pr}_2\circ d\Psi_v=\hbox{\rm natural projection}.\cr}\right.$$
\smallskip
2) We define a chart $\Psi^*:V^*\to T_L{\cal L}^*_v\times N_L{\cal L}^*_v$ such that 
$\Psi^{*-1}(\alpha,0)$ passes through
$v$ and $\Psi^{*-1}(0,\beta)$ is ``horizontal''. More precisely, $\Psi^{*-1}(\alpha,\beta)$
 is given in the chart $\Psi$
by an equation 
$$y=f_{\alpha,\beta}(x)\,\,,\,\,x\in T_vL\,\,,\,\,y\in N_vL,$$ 
such that
$$\left\{\eqalign{f_{\alpha,\beta}(0)&=\nu(\beta)\cr
f_{\alpha,0}(0)&=0\cr
{\partial f_{\alpha,0}\over\partial x}(0)&=i(\alpha).\cr}\right.$$
In the last equation, $\alpha\in T_L{\cal L}^*_v=P^*$ is interpreted as an element of $T_PE_v$, so that $i(\alpha)\in {\rm Hom}(T_vL,N_vL)$. 
\smallskip
3) Finally, let $P^*\in Gr_2T_LV^*$ close to $T_LV^*$, we define 
$\chi(P^*)\in{\rm Hom}(T_L{\cal L}^*_v,N_L{\cal L}^*_v)$ as the unique
$h$
such that $$d\Psi^*_L(P^*)={\rm graph}(h).$$

\medskip
\noindent{\bf End of the proof.}
\hskip3mm
Let $w=\Psi^{-1}(x)$ be an element of $L$ close to $v$. Then $\Psi^{*-1}(\alpha,\beta)\in V^*_w$ 
if and only if $f_{\alpha,\beta}(x)=0$, 
thus $T_LV^*_w$ is given by $$\{(\delta\alpha,\delta\beta)\mid 
{\partial f_{\alpha,0}\over\partial\alpha}(x).\delta\alpha
+{\partial f_{0,\beta}\over\partial\beta}(x).\delta\beta=0\},$$
ie $${\partial f_{\alpha,0}\over\partial\alpha}(x).\delta\alpha+\nu(\delta\beta)=0.$$
In other words
$$\chi( T_LV^*_w)=-\nu^{-1}({\partial f_{\alpha,0}\over\partial\alpha}(x))$$
Thus, the tangent space of $E^*_L$ at $T_L{\cal L}^*_v$
is identified with the image of the morphism
$$i^*=\nu^{-1}\circ{\partial^2 f_{\alpha,0}\over\partial x\partial\alpha}(0):T_vL\to {\rm Hom}(T_L{\cal L}^*_v,N_L{\cal L}^*_v).$$
\smallskip
\noindent Since ${\partial f\over\partial x}(\alpha,0,0)=i(\alpha)$, one has
$$i^*(\xi) 
(\alpha)=i(\alpha) 
(\xi)\,\,,\,\,(\xi,\alpha)\in T_vL\times T_L{\cal L}^*_v.$$
Since $E_v$ is elliptic, $i(\alpha)$ is invertible and orientation-preserving if $\alpha\ne 0$. 
Thus one can identify the oriented planes $T_L{\cal L}^*_v$, $T_vL$ and $N_L{\cal L}^*_v$ 
with $\C$ so that $i(\alpha)$ is the 
multiplication by $\alpha$. Then $i^*(\xi)$ is the multiplication by $\xi$, thus
it  is invertible and orientation-preserving 
if $\xi\ne0$, which means that $E^*_L$ is elliptic.
\bigskip
\noindent{\bf 4.C. Proof that $V^*$ is oriented diffeomorphic to $\C\P^2$.} \hskip3mm 
One could prove it similarly to the proof for $V$. The simplest proof however is to remark that
the space of elliptic structures on $V=\C\P^2$ which are tamed by $\omega_0$ is contractible. 
For each $E$ in this space, we obtain an
oriented manifold $V^*_E$ which varies smoothly with $E$, thus keeps the same oriented diffeomorphism type. 
Since for $E$ associated to $J_0$ one has $V^*_E=\C\P^{2*}$ the standard dual projective plane, this proves the result.

\bigskip
\noindent{\bf 4.D. Tameness of $E^*$ and identification $(V^*)^*=V$.}\hskip3mm For each $v\in V$, the surface 
${\cal L}^*_v\subset V^*$ is an $E^*$-curve of degree $1$, ie an $E^*$-line. Moreover, for two distinct points $L,L'\in V^*$
there exists a unique $v\in L\cap L'$, equivalently a unique ${\cal L}^*_v$ containing $L$ and $L'$: this means that the 
$E^*$-lines are precisely the ${\cal L}^*_v$, and thus that $E^*$ is tame and $V^{**}=V$. The equivalence
$(v\in L\Leftrightarrow L\in{\cal L}^*_v)$ implies that $E^{**}$
 is identified to $E$.
\bigskip
\noindent{\bf 4.E. Dual curves.}\hskip3mm Let ${\cal J}$ be the restriction to $TI$ of $(J,J^*)$, where
$J$ and $J^*$ are the twisted almost complex structures associated to $E$ and $E^*$: it is an almost complex 
structure, whose images by $\gamma$ and $\gamma^*$ (notations of 4.A) are  
 $\widetilde J$ and  $\widetilde J^*$, the complex structures on $\Theta$ and $\Theta^*$ associated to $E$ and $E^*$. 
Thus the map $\phi:(E,\Theta)\to(E^*,\Theta^*)$ is a $(\widetilde J,\widetilde J^*)$-biholomorphism.
\smallskip
Now let $C=f(S)\subset V$ be an irreducible $E$-curve (or an irreducible germ) not contained in an $E$-line. Let $\gamma:S\to E$
be the Gauss map, which is $\widetilde J$-holomorphic. Then
$\gamma^*=\phi\circ\gamma:S\to E^*$  is $\widetilde J$-holomorphic and not locally constant, thus it is the Gauss map
of an $E^*$-map $f^*:S\to V^*$. By definition, $C^*=f^*(S)$ is the {\it dual} curve of $C$: it is again an irreducible
$E^*$-curve (or germ), not contained in an $E^*$-line, and of course one has $C^{**}=C$.
\bigskip
\bigskip
\noindent{\bf 5. Inexistence of a natural almost complex structure on $V^*$}
\medskip
Here we construct a tame almost complex structure $J$ on $V=\C\P^2$ 
such that the twisted almost complex structure $J^*$ on $V^*$ is non linear.

\medskip
We impose on $J$ the following properties:
\smallskip
- it is standard outside $$U_0=(\Delta(2)\setminus\Delta(1))\times\Delta(2)\subset\C^2\subset\C\P^2$$

- it is $\omega_0$-positive
\smallskip
- for $\alpha\in\C$ small enough, the $J$-line $L(\alpha,\alpha)$ passing through the point $(0,\alpha)$ with the slope
$\alpha$ has an intersection with $U_0$ given by the equation
$$y=f_{\alpha,\alpha}(x)=\alpha+\alpha x+{1\over5}\rho(x).\overline\alpha x^2,$$
where $\rho:\C\to[0,1]$ takes the value $1$ on $U_1=\Delta(1)\times\Delta(1)$ and $0$ outside $U_0$. 
Note that the factor ${1\over5}$ guarantees that $\alpha\mapsto f_{\alpha,\alpha}(x)$ is an embedding for $|x|<2$
near $0$.
\smallskip
One can find such a $J$ under the form
$$J(x,y)=\pmatrix{i&0\cr
b(x,y)\sigma&i}$$
where $b(x,y)\in\C$ and $\sigma$ is the complex conjugation. Then $L(\alpha,\alpha)\cap U_0$ is $J$-holomorphic
if and only if
$$b(x,f_{\alpha,\alpha}(x))={\partial f_{\alpha,\alpha}\over\partial\overline x}(x).$$
Since $\alpha\mapsto f_{\alpha,\alpha}(x)$ is an embedding near $0$ for $|x|<2$ and the second member
vanishes for $x|$ close to $2$, one can find a smooth $b(x,y)$ with support in $U_0$, satisfying the above equality for
$|x|<2$ and $\alpha$ small enough.
\medskip
We now prove that $J^*_L$ is not linear. Note that on $U_1$, we have $J=J_0$ and $L(\alpha,\alpha)\cap U_1$ is given
by
$$y=f_{\alpha,\alpha}(x)=\alpha+\alpha x+{1\over5}\overline\alpha x^2.$$
Let $L$ be the  $J$-line $L(0,0)$, which is the $x$-axis. Recall that for each $v\in L$ the subspace $T_L{\cal L}^*_v\subset
T_LV^*$ is preserved by $J^*-L$, which is linear on it.
The global linearity of $J^*_L$ is equivalent to the following:
$$\matrix{(\forall\xi,\eta\in T_LV^*)&\xi+\eta\in T_L{\cal L}^*_v\Rightarrow J^*_L(\xi)+J^*_L(\eta)\in T_L{\cal L}^*_v.}$$
\smallskip
Consider on $L$ the points $v_0=0$ and $v_1=\infty$. Then we have a direct sum
$T_LV^*=T_L{\cal L}^*_0\oplus T_LL^*_\infty$. Then fix $\alpha\ne0$ and consider the path $t\in[0,1]\mapsto \gamma(t)=L(t\alpha,t\alpha)\in V^*$,
and write its derivative at $t=0$ as $$\dot\gamma=\xi+\eta\,\,,\,\,\xi\in T_L{\cal L}^*_0\,\,,\,\,T_L{\cal L}^*_\infty.$$
It belongs to $T_L{\cal L}^*_v$ where $v=\lim_{t\to0}(\gamma(t)\cap L)$. Identifying $L$ with $\C\P^1$, 
this means that 
$v$ is the solution of the equation
$$\alpha+\alpha v+{1\over5}\overline\alpha v^2=0.$$
If we change $\alpha$ to $i\alpha$, $\xi$ and $\eta$ are changed to $J^*(\xi)$ and $J^*(\eta)$ (essentially since $J_0$ is standard
near $0$ and $\infty$, thus $J^*_L(\xi)+J^*_L(\eta)\in T_L{\cal L}^*_w$ where $w$ is the solution of 
$$\alpha+\alpha w-{1\over5}\overline\alpha w^2=0.$$
Thus $w\ne v$, which means that $J^*_L$ is not linear.

\smallskip

\medskip
\noindent{\bf Conjecture.}\hskip3mm {\it Let $(V,J)$ be a tame almost complex projective plane. 
Assume that  the elliptic structure $E^*$ on $V^*$ comes from an almost complex structure $J^*$
 on $V^*$.
Then $J$ is integrable, thus $(V,J)$ is biholomorphic to $(\C\P^2,J_0)$.
 \smallskip More precisely: if $J^*_L$ is linear, then 
the Nijenhuis torsion of $J$ vanishes on $L$.}
\bigskip
\bigskip

\noindent{\bf 6. Pl\"ucker formulas for $E$-curves}
\medskip
We  follow the classical
topological method in algebraic geometry, cf. for instance [GH] p.279.
\medskip
Let $C=f(S)\subset V$ be an irreducible $E$-curve, not contained in an $E$-line, and
let $C^*\subset V^*$ be its dual. We compute the degree $d^*$ of $C^*$, which is the number of intersection points of $C^*$ with an $E^*$-line, ie the number of points
of $C$ such that the tangent line $L_vC$ contains $v$. This number is to be interpreted algebraically, but for a generic
$v$ it is equal to the set-theoretic number. 
\smallskip
Let $L$ be an $E$-line disjoint from $v$, then the central projection $V\setminus\{v\}\to L$ along $E$-lines 
through $v$
induces an ``almost holomorphic''
branched covering
$C\to L$ of degree $d$, in the sense that each singularity has a model $z\to z^k$: this is a consequence of the positivity 
of intersections. Let $S$ be the normalization of $C$, then the number of branch points of 
the induced covering $S\to L$ is $d^*+\kappa$ where $\kappa$ is the algebraic number of cusps, ie
the algebraic number of zeros of $df$ if $f$ is a parametrization of $C$.
Thus we have the
Hurwitz formula $2-2g=2d-(d^*+\kappa)$, where $g$ is the genus of $C$, ie
$$d^*=2d+2g-2-\kappa.$$
In particular, if $C$ has only $\delta$ nodes and $\kappa$ cusps, we have $2g-2=d(d-3)-2\delta-2\kappa$ thus we get the first
Pl\"ucker formula
$$d^*=d(d-1)-2\delta-3\kappa.$$
As in the classical case, the other Pl\"ucker formulas follow from this and the genus formula, with the fact that 
an ordinary bitangent (resp. flex) of $C$ corresponds to a node (resp. cusp) of $C^*$.
\medskip
This  implies restrictions on the possible sets of singularities going beyond the genus formula.
 For instance, if $C$ has only  nodes and cusps, then another form of Pl\"ucker formula is
$$\kappa=2g-2+2d-d^*.$$
If $d=5$ and $g=0$ we get $\kappa=8-d^*$, and since $d^*\ge3$ we have $\kappa\le5$: 
not all $6$ nodes of a generic rational curve can be transformed to cusps.

\medskip In general, if $C$ is rational with only nodes and cusps, we get $\kappa=2d-2-d^*<3d$, which implies that
the space of rational $J$-curves is, at the point $C$, a smooth manifold of the expected dimension 
(equal to $d(d+3)$ over $\R$): this follows from the generalization of the 
automatic genericity proved in [B]. This could be interesting for the isotopy problem of symplectic surfaces 
[Sik2].

\vskip10mm

\centerline{\bf References}
\medskip
\noindent[AL]  M. Audin and J. Lafontaine eds, {\it Holomorphic curves
in symplectic geometry}, Progress in Math. 117, Birkhauser, 1994.
\medskip
\noindent[Aur] D. Auroux,  {\it Symplectic 4-manifolds as branched coverings of $\C\P^2$},
 Invent. Math. {\bf 139} (2000), 551-602.
\medskip
\noindent{[B]} J.-F. Barraud, {\it Nodal symplectic spheres in $\C\P^2$ with positive
self-intersection}, 
Intern. Math. Res. Not. {\bf9} (1999), 495-508.
\medskip
\noindent [G]  M. Gromov, {\it Pseudo holomorphic curves in symplectic
manifolds},
Invent. Math. {\bf82} (1985), 307-347.
\medskip
\noindent [GH] P. Griffiths 	and J. Harris, {\it
Principles of algebraic geometry}, John Wiley $\&$ Sons, New York, 1978.
\medskip
\noindent [HLS]  H. Hofer, V. Lizan, J.-C. Sikorav, {\it On genericity for
holomorphic
curves in $4$-dimensional almost-complex manifolds}, J. Geom. Anal. {\bf7} (1998), 149-159.
\medskip
\noindent
[L], F. Labourie, {\it 
Probl\`emes de Monge-Amp\`ere, courbes holomorphes et laminations},
Geom. Funct. Anal. 7 (1997), 496-534. 
\medskip
\noindent [McD1] D. McDuff, {\it Rational and ruled symplectic $4$-manifolds}, 
J. Amer. Math. Soc. {\bf 3} (1990), 679-712.
\medskip
\noindent [McD2] D. McDuff, {\it The local behaviour of holomorphic curves in almost complex $4$-manifolds}, 
J. Diff. Geom. {\bf 34} (1991), 143-164.
\medskip
\noindent [MW] M. Micallef and B. White, {\it The structure of branched points in minimal
surfaces and in pseudoholomorphic curves}, Ann. of Math. {\bf 139} (1994), 35-85.
\medskip
\noindent [Sik1], J.-C. Sikorav, {\it Singularities of $J$-holomorphic
curves}, Math. Z. {\bf226} (1997), 359-373.
\medskip
\noindent [Sik2], J.-C. Sikorav, {\it The gluing construction for normally generic $J$-holomorphic curves}, 
preprint Ecole Normale Sup\'erieure de Lyon, UMPA, {\bf 264} (2000).
\medskip
\noindent [T] C. Taubes, {\it ${\rm SW}\Rightarrow{\rm Gr}$: from the Seiberg-Witten equations to 
pseudo-holomorphic curves}, 
J. Amer. Math. Soc. 9 (1996), no. 3, 845--918. 
\medskip
\noindent [V], I.N. Vekua, {\it Generalized analytic functions},
Pergamon Press, 1962.

\vskip10mm
Jean-Claude Sikorav

Ecole Normale Sup\'erieure de Lyon,

Unit\'e de Math\'ematiques Pures et Appliqu\'ees

46, all\'ee d'Italie

69364 Lyon cedex 07, FRANCE

sikorav@umpa.ens-lyon.fr

\end